\documentclass[12pt]{article}
\usepackage{amsmath,amsthm,color}
\allowdisplaybreaks

\theoremstyle{plain}

\begin{document}
\baselineskip5pt
\title{\Large{A Short Proof that Lebesgue Outer Measure of an Interval is Its Length}}
\author{Jitender Singh\\ \small Department of Mathematics, Guru Nanak Dev University Amritsar, India\\ \small \color{magenta}\tt sonumaths@gmail.com; jitender.math@gndu.ac.in}
\date{\empty}
\maketitle
\baselineskip12pt
\noindent
The Lebesgue outer measure $m^*(E)$ of a subset $E$ of real line is defined as  $m^*(E):=\inf\{\sum_{k=1}^\infty \ell(I_k)~|~E\subseteq \cup_{k=1}^\infty I_k\}$, where each $I_k$ is an open interval and $\ell(I_k)$ is its length. Establishing one of the inequalities in the standard  proof of the fact in the title above turns out to be tedious in \cite[p.\ 31]{p1}. Using the connectedness of the interval shortens the proof as follows.
\begin{proof}
Given two real numbers $a$ and $b$ with $a<b$, it is enough to prove that $m^*([a,b])=b-a$.
Clearly, $m^*([a,b])\leq b-a$. Now let  $[a,b]\subset \cup_{k=1}^n I_k$ for some positive integer $n$, which is always possible since $[a,b]$ is compact. Without loss of generality, assume that the set $[a,b]\cap I_k$ is nonempty for each $k$. Observe that  the set $\cup_{k=1}^n I_k$ is connected. (Otherwise, if $(P,Q)$ is its separation, then for each $k$, by connectedness of $I_k$, either $I_k\subset P$ or $I_k\subset Q$. Thus each of $P$ and $Q$ is equal to union of sets from the list $\{I_1,\ldots,I_n\}$. So the pair $(P\cap[a,b],Q\cap [a,b])$ determines a separation of $[a,b]$, which contradicts connectedness of $[a,b]$.) So $\cup_{k=1}^n I_k$ is an open interval containing $[a,b]$. Thus, $b-a\leq \ell(\cup_{k=1}^n I_k)\leq \sum_{k=1}^n\ell(I_k)$, where the last inequality holds since some intervals overlap$^\dagger$. Hence, $b-a\leq m^*([a,b])$.
\end{proof}
$^\dagger$The inequality $\ell(\cup_{k=1}^n I_k)\leq \sum_{k=1}^n\ell(I_k)$ can be justified as follows. Given a bounded interval $I$, by definition, $\ell(I)=\sup I- \inf I$. Observe that if two distinct bounded intervals $I_1$ and $I_2$ overlap, then $\ell(I_1\cup I_2)\leq \ell(I_1)+\ell(I_2)$.

Given $n>1$ bounded open intervals $I_1,\ldots I_n$ with their union being connected implies that for each $k=1,\ldots,n-1$  we may choose $I_k$ after re-indexing  these intervals, such that $\cup_{j=1}^kI_j$ is connected and it overlaps at least one interval among rest of the $(n-k)$ intervals, which we denote by $I_{k+1}$. So, we  have $\ell (\cup_{k=1}^{^n}I_k)=\ell(\cup_{k=1}^{n-1}I_k\cup I_n)\leq \ell(\cup_{k=1}^{n-1}I_k)+\ell(I_n)\leq \ell(\cup_{k=1}^{n-2}I_k)+ \ell(I_{n-1})+\ell(I_n)\leq\ldots$ $\leq \ell(I_1)+\ldots+\ell(I_n)$.


\begin{thebibliography}{1}
\bibitem{p1} H.\ L.\ Royden, P.\ M.\ Fitzpatrick, \textit{Real Analysis}. Fourth ed. Pearson, Boston, 2010.
\end{thebibliography}
\end{document}